\documentclass[reqno,12pt]{amsart}
\usepackage{graphics}
\usepackage{a4,amsmath,epsfig}

\def\goth#1{{\mathfrak #1}}

\def\0{\goth 0}
\def\gA{\goth A}
\def\gB{\goth B}

\def\Id{\rm Id}
\begin{document}

%%%%%%%%%%%%%%%%%%%%%%%%%%%%%%%%%%%%%%%%%%%%%%%%%%%%%%%%%%%%%%%%%%%%%%%%%%%%%%%
\vspace*{0.2cm}

\begin{center}
{\Large \bf On principal fibrations associated with one algebra}\\[1.5 cm]

{\bf\large Maria Trnkov\'a$^{1,2}$}\\[0.4cm] email: {\tt
Maria-Trnkova@seznam.cz}\\[0.3cm]
{\it$^1$Department of Algebra and Geometry, Faculty of Science,}\\
{\it Palacky University in Olomouc, Tomkova 10, Olomouc, Czech
Republic}\\[0.1cm]
{\it$^2$Department of Geometry, NIMM,}\\
{\it Kazan State University, prof. Nuzhina 1/37, Kazan, Russia}
\\[0.5 cm]

\end{center}

To the memory of B.N. Shapukov

\vspace*{1.0cm}

\textbf{Abstract} \emph{In this paper we study two types of
fibrations associated with a 3-dimensional unital associative
irreducible algebra and their basic properties. We investigate
trivial principal fibrations of degenerate semi-Euclidean sphere and
their semi-conformal and projective models. We use Norden
normalization method for constructing second model.}

\section{Principal fibrations of $G$ by its subgroups}
Let us denote by $\gA$ a unital associative $n$-dimensional algebra
with multiplication $xy$, by $G\subset\gA$ the set of invertible
elements. It is well known that $G$ is a Lie group with the same
multiplication rule. Let ${\gB}$ be a unital subalgebra of algebra
${\gA}$ and $H\subset G$ the set of invertible elements. So $H$ is a
Lie subgroup of group $G$. We consider the factor-space $G/H$ of
right cosets. A fibration $(G, \pi, M=G/H)$ is a principal fibration
with the structure group $H$, where $\pi$ is a canonical projection
\cite{Shap}. Here we have a fibration:
$$H \rightarrow G \rightarrow G/H.$$

We consider all 3-dimensional unital associative irreducible
algebras up to isomorphism. It is known that there exist only tree
types of such algebras \cite{Stud}, \cite{Vish}. In special choice
of basic units they have the following multiplication rules:
$$
{\rm I})\begin{tabular}{c|c|c|c|}
 & $1$ & $e_{1}$ & $e_{2}$  \\
\hline
$1$ & $1$ & $e_{1}$ & $e_{2}$  \\
\hline
$e_{1}$ & $e_{1}$ & $e_{2}$ &  $0$ \\
\hline
$e_{2}$ & $e_{2}$ &  $0$ & $0$ \\
\hline
\end{tabular}
\qquad {\rm II})\ \begin{tabular}{c|c|c|c|}
 & $1$ & $e_{1}$ & $e_{2}$  \\
\hline
$1$ & $1$ & $e_{1}$ & $e_{2}$  \\
\hline
$e_{1}$ & $e_{1}$ & $1$ &  $e_{2}$ \\
\hline
$e_{2}$ & $e_{2}$ &  $-e_{2}$ & $0$ \\
\hline
\end{tabular}
\qquad {\rm III})\ \begin{tabular}{c|c|c|c|}
 & $1$ & $e_{1}$ & $e_{2}$  \\
\hline
$1$ & $1$ & $e_{1}$ & $e_{2}$  \\
\hline
$e_{1}$ & $e_{1}$ & $0$ &  $0$ \\
\hline
$e_{2}$ & $e_{2}$ &  $0$ & $0$ \\
\hline
\end{tabular}
.$$

For these algebras N.Belova found all subalgebras with unit and
principal fibrations \cite{Bel}. Only in algebras of type $\rm II$
and $\rm III$ exist the assosiated conjugation
$x=x_{0}+x^{i}e_{i}\to \overline{x}=x_{0}-x^{i}e_{i}$ which has the
property $\overline{xy}=\overline{y}\,\overline{x}$.

We consider the bilinear form
\begin{equation}\label{1}
(x,y)=\frac{1}{2}(x\overline{y}+y\overline{x}).
\end{equation}
For algebras of type $\rm II$ and $\rm III$  this form takes the
real values and it determines a degenerate scalar product on them:
$$\rm II) (x,y)=x_{0}y_{0}-x_{1}y_{1},\quad \rm III) (x,y)=x_{0}y_{0}.$$ So
these algebras have a structure of semi-Euclidean vector spaces with
rank 2 and 1, respectively. In our discussion we concentrate on the
algebra of type $\rm II$, because it is less degenerate than type
$\rm III$.

This is non-Abelian algebra. Multiplication rule for algebra
elements is:
\begin{equation}\label{xy}
xy=(x_{0}+x_{1}e_{1}+x_{2}e_{2})(y_{0}+y_{1}e_{1}+y_{2}e_{2})=
\end{equation}
$$x_{0}y_{0}+x_{1}y_{1}+(x_{0}y_{1}+x_{1}y_{0})e_{1}+(x_{2}(y_{0}-y_{1})+(x_{0}+x_{1})y_{2})e_{2}.$$
An inverse element of $x$ is:
\begin{equation}\label{rever}
x^{-1}=\frac{x_{0}-x_{1}e_{1}-x_{2}e_{2}}{(x_{0})^{2}-(x_{1})^{2}}.
\end{equation}
The set of invertible elements $G=\{x\in {\gA}\mid
(x_{0})^{2}-(x_{1})^{2}\ne 0\}$ is a non-Abelian Lie group. Its
underlying manifold is $\mathbb{R}^{3}$ without two transversal
2-planes, hence it consists from 4 connected components.

\textbf{Proposition}(N.Belova) \emph{Any 2-dimensional subalgebra of
algebra of type $\rm II$ with an algebra unit is isomorphic to
double or dual numbers algebras.}

We consider a subalgebra $R(e_1)$ with basis $\{1, e_1\}$, it is an
algebra of double numbers, and a subalgebra $R(e_2)$ with basis
$\{1, e_2\}$, it is an algebra of dual numbers. The set of their
invertible elements $H_{1}=\{x_{0}+x_{1}e_{1}\in R(e_{1})\mid
x^{2}_{0}-x^{2}_{1}\ne 0\}$ and $H_{2}=\{x_{0}+x_{2}e_{2}\in
R(e_{2})\mid x_{0}\ne 0\}$ are Lie subgroups of the Lie group $G$.
First we take the space of right cosets by $R(e_1)$ subalgebra and
fibration by it. Then we have the following proposition:

\textbf{Proposition} (N.Belova) \emph{Fibration $(G, \pi, M=G/H_1)$
determined by formula
\begin{equation}\label{pi1} \pi(x)=\frac{x_{2}}{x_{0}-x_{1}}
\end{equation}
is a trivial principal fibration over the real line $\mathbb{R}$.
The typical fiber, it is a plane without two transversal lines. The
structure group is $H_1$.}

Therefore, the manifold of the group $G$ is diffeomorphic to direct
sum $\mathbb{R}\times H_{1}$. The equation of fibers is:
\begin{equation}\label{fiber1} u(x_{0}-x_{1})-x_{2}=0, \; u\in R.
\end{equation}
This is 1-parametric family of planes with common axis:
$x_{0}-x_{1}=0,\,x_{2}=0$. Of course it is necessary to remove the
intersection points of these planes with planes: $x_{0}\pm x_{1}=0$.

We consider the left multiplications $x'=ax$ on the group $G$. They
form the 3-parametric Lie group of linear transformations:
\begin{equation}\label{L1}
L(a)=\left(\begin{array}{crc} a_{0}&a_{1}&0\\
a_{1}&a_{0}&0\\
a_{2}&-a_{2}&a_{0}+a_{1}\\
\end{array} \right),
\end{equation}
where determinant
${\det}L(a)=(a^{2}_{0}-a^{2}_{1})(a_{0}+a_{1})\neq0$. The
multiplications preserve fibers if and only if when $a\in H_{1}$.
This group consists of four connected components.

The right multiplications $x'=xb$ on the group $G$ forms the
3-parametric Lie group of linear transformations:
\begin{equation}\label{R1}
R(b)=\left(\begin{array}{crc} b_{0}&b_{1}&0\\
b_{1}&b_{0}&0\\
b_{2}&b_{2}&b_{0}-b_{1}\\
\end{array} \right)
.\end{equation}
 They preserve the fibration and induce the 2-parametric group of affine
transformations on the base:
$$
u'=\alpha u+\beta,
$$
here $\alpha=\frac{b_{0}-b_{1}}{b_{0}+b_{1}}\ne 0$, $\beta=\pi(b)$.

The right multiplications include involutions of hyperbolic type
with matrices:
$$
R(b_2)=\left(\begin{array}{crc} 0&1&0\\
1&0&0\\
b_{2}&b_{2}&-1\\
\end{array} \right),\quad
 R^{2}={\Id}.
$$
They induce the 1-parametric family involutions on the base:
$$
u'=-u+b_2.
$$

Let us consider rotations and anti-rotations of this semi-Euclidean
space. Multiplications $x'=ax$ and $x'=xa$, where $|a|^2=\pm1$, are
rotations and anti-rotations because of the relations:
$$|x'|^2=|a|^2|x|^2=\pm|x|^2.$$ Each element $a\in\gA$,
$|a|^2=\pm1$, can be represented as:
$$a=\cosh\varphi+a_0\sinh\varphi \quad \mbox{ or } \quad
a=\sinh\varphi+a_0\cosh\varphi,$$ where
$a_0=-\overline{a}_0,\;|a_0|^2=-1.$ The bilinear form (\ref{1}) in
the algebra $\gA$ takes real values, therefore it is possible to
present it as:
$(x,y)=\frac{1}{2}(x\overline{y}+y\overline{x})=\frac{1}{2}(\overline{x}y+\overline{y}x).$
Consequently, the hyperbolic cosines or sines of angle between $x$
and $x'$ is equal:
$$
\frac{(x,ax)}{|x||ax|}=\frac{1/2(x\overline{ax}+ax\overline{x})}{|x|^2}=
\frac{1/2(x\overline{x}\,\overline{a}+ax\overline{x})}{|x|^2}=
\frac{1}{2}(\overline{a}+a),
$$
$$
\frac{(x,xa)}{|x||xa|}=\frac{1/2(\overline{x}xa+\overline{xa}x)}{|x|^2}=
\frac{1/2(\overline{x}xa+\overline{a}\,\overline{x}x)}{|x|^2}=
\frac{1}{2}(a+\overline{a}).
$$
In the last two equations we get $\cosh\varphi$ if $|a|^2=1$ and we
get $\sinh\varphi$ if $|a|^2=-1$. It means that this angle does not
depend on $x$.

Transformations
\begin{equation}\label{rot1}
x'=axb,
\end{equation}
where $|a|^2=\pm1,\; |b|^2=\pm1$, are compositions of rotations or
anti-rotations $x'=ax$ and $x'=xb$. So, there are proper rotations
or anti-rotations.

Transformations
\begin{equation}\label{rot2}
x'=a\overline{x}b
\end{equation}
are compositions of (\ref{rot1}) and reflection $x'=\overline{x}$.
So, there are improper rotations or anti-rotations.

\textbf{Lemma} \emph{Any proper or improper rotation and
anti-rotation of this semi-Euclidean space can be represented by
(\ref{rot1}) or (\ref{rot2}).}

\textbf{Proof} These rotations and anti-rotations are compositions
of odd and even numbers of reflections of planes passing through the
origin. To each plane corresponds its orthonormal vector $n$. If
vectors $x_1$ and $n$ are collinear, then
$\overline{x}_1n=\overline{n}x_1$ and
$x'_1=-n\overline{x}_1n=-n\overline{n}x_1=-x_1$. If vectors $x_2$
and $n$ are orthogonal, then $\overline{x}_2n+\overline{n}x_2=0$ and
$x'_2=-n\overline{x}_2n=n\overline{n}x_2=x_2$. On the other hand any
vector $x$ can be represented as a sum of vectors $x_1$ and $x_2$.
It means, that a reflection of plane is: $x'=-n\overline{x}n$. So,
composition of odd and even numbers of planes reflections are
transformation (\ref{rot1}) or (\ref{rot2}). $\Box$

Let us introduce adapted coordinates $(u,\lambda,\varphi)$ to
fibration in semi-Euclidean space, here $u$ is a basic coordinate,
$\lambda,\varphi$ are fiber coordinates. If $|x|^{2}>0$, we denote
$\lambda =\pm\sqrt{x_{0}^{2}-x_{1}^{2}}\ne 0$, the sign of $\lambda$
is equal to the sign of $x_0$. The adapted coordinates of fibration
in this case are:
\begin{equation}\label{adapt>0}
x_{0}=\lambda \cosh \varphi ,\quad x_{1}=\lambda \sinh \varphi ,
\quad x_{2}=u\lambda \exp \varphi ,
\end{equation}
where $\lambda \in\mathbb{R}_{0},\quad u,\varphi\in\mathbb{R}.$

If $|x|^{2}<0$, then we write $\lambda
=\pm\sqrt{x_{1}^{2}-x_{0}^{2}}$, the sign of $\lambda$ is equal to
the sign of $x_1$:
\begin{equation}\label{adapt<0}
x_{0}=\lambda \sinh \varphi ,\quad x_{1}=\lambda\cosh\varphi ,\quad
x_{2}=u\lambda \exp \varphi .
\end{equation}

The structure group acts as follows:
\begin{equation}\label{str-gr1}u'=u, \quad
  \lambda'=\lambda \rho, \quad
  \varphi'=\varphi+\psi,
\end{equation}
where the element $a(0, \rho, \psi)$ of structure group acts on the
element $x(u, \lambda, \varphi)\in G$. This group consists of 4
connected components.

\section{The principal subfibration of the fibration $(G,\pi,M=G/H_1)$}
As we already said, the scalar product in the algebra $\gA$ of type
$\rm II$ is: $(x,y)=x_{0}y_{0}-x_{1}y_{1}$. So, the algebra $\gA$ is
a 3-dimensional semi-Euclidean space with rank 2. We call
\emph{semi-Euclidean sphere with unit radius} the set of all
elements of algebra $\gA$ whose square is equal to one, $$
S^{2}(1)=\{x\in\gA\mid x_{0}^{2}-x_{1}^{2}=1\}.
$$ It looks like a hyperbolic cylinder. The set of elements with imaginary
unit module $|x|^{2}=-1$ we call \emph{semi-Euclidean sphere with
imaginary unit radius} $S^{2}(-1)$. One of these spheres can be
obtained up from another one by rotation.

We consider the subfibration of the fibration $(G,\pi,M=G/H_1)$ to
semi-Euclidean sphere $S^{2}(1)$, i.e. the fibration  $\pi:
S^{2}(1)\to M$. The fibers of new fibration are intersections of
$S^2(1)$ and planes (\ref{fiber1}). The restriction of the group of
double numbers $H_1$ to $S^{2}(1)$ is a Lie subgroup $S_1$ of double
numbers with unit module
$$ S_{1}=\{a_{0}+a_{1}e_{1}\in H_{1}\mid a_{0}^{2}-a_{1}^{2}=1\}\,.
$$ This group consists of two connected components.

\textbf{Proposition} \emph{The fibration $(S^{2}(1),\pi,M)$ is
principal fibration of the group $S^{2}(1)$ by the Lie subgroup
$S_1$ to right cosets.}

\textbf{Proof.} Let $x$ and $y$ be two sphere points from one fiber.
They are also from the same fiber of the principal fibration
$(G,\pi,M)$. So, there exists a unique element $a\in H_{1}$ such
that $y=ax$. Then $|y|^{2}=|a|^{2}|x|^{2}$ and hence $|a|^{2}=1$. It
means that $a\in S_{1}$. $\Box$

We define coordinates adapted to the fibration on semi-Euclidean
sphere $S^{2}(1)$. If $x\in S^{2}(1)$ then by (\ref{adapt>0}) we get
$\lambda =\varepsilon,\, \varepsilon =\pm 1$. The parametric
equation of semi-Euclidean sphere in the adapted coordinates
$(u,\varphi)$ is:
\begin{equation}\label{adapt1}
{\textbf{r}}(u,\varphi)=\varepsilon(\cosh \varphi, \sinh
\varphi,u\exp\varphi),
\end{equation}
where $u$ is a basis coordinate, $\varphi$ is a fiber coordinate.
Different values of $\varepsilon$ correspond to different connected
components of semi-Euclidean sphere $S^{2}(1)$.

Let us define the action of the structure group $S_{1}$ on
semi-Euclidean sphere. By (\ref{str-gr1}) and the adapted
coordinates of elements $a(0, \varepsilon_{1}, \psi),\, x(u,
\varepsilon, \varphi)\in S^{2}(1)$ we get:
$$u'=u, \quad
  \varepsilon'=\varepsilon \varepsilon_{1}, \qquad
  \varphi'=\varphi+\psi. \qquad
$$
This group also consists of two connected components.

The metric tensor for semi-Euclidean sphere has the matrix
representation:
$$
(g_{ij})=\left(\begin{array}{cc}
0 & 0 \\
0 & -1
\end{array}
\right),\quad \mathrm{rank} (g_{ij})=1
$$
So the linear element of metric is:
\begin{equation}\label{metric1}
ds^{2}_{1}=-d\varphi^{2}.
\end{equation}

\section{Semi-conformal model of the sphere fibration $(S^{2}(1),\pi,\mathbb{R})$}

We consider the semi-conformal model of the fibration
$(S^{2}(1),\pi,\mathbb{R})$. We project stereographicaly the sphere
$S^{2}(1)$ from the point $N(1, 0 ,0)\in S^{2}(1)$ (pole) to the
equatorial plane $\mathbb{R}^{2}$ with equation $x_{0}=0.$ We write
the stereographic map in the coordinate form. For this we find  the
equation of line which pass through the pole and an arbitrary point
of the sphere $X(x_{0},x_{1},x_{2}):\quad
\widetilde{x}_{0}=1+(x_{0}-1)t,\,\widetilde{x}_{1}=x_{1}t,\,\widetilde{x}_{2}=x_{2}t$.
 When this line crosses the equatorial plane, we get the formulas of
stereographic map $f:S^{2}(1)\rightarrow \mathbb{R}^{2}$ where
$x_{0}\neq 1$:
\begin{equation}
\label{stereo1}x=\frac{x_{1}}{1-x_{0}},\quad
y=\frac{x_{2}}{1-x_{0}}\,,
\end{equation}
here $(x,y)\in \mathbb{R}^{2}, \quad (x_{0},x_{1},x_{2})$ are
coordinates on $S^{2}(1)$. An inverse map
$f^{-1}:\mathbb{R}^{2}\rightarrow S^{2}(1)$ when $x\neq \pm 1$ is:
\begin{equation}
\label{rev-stereo1}x_{0}=-\frac{1+x^{2}}{1-x^{2}},\quad
x_{1}=\frac{2x}{1-x^{2}},\quad x_{2}=\frac{2y}{1-x^{2}}\,.
\end{equation}
If we put formulas (\ref{stereo1}) into (\ref{adapt1}) then we
obtain the relations between coordinates $x, y$ and adapted
coordinates $u,\varphi$ which are on semi-Euclidean sphere:
$$f:\quad x=\frac{\sinh
\varphi}{\varepsilon-\cosh \varphi},\quad
y=\frac{u\exp\varphi}{\varepsilon-\cosh \varphi}\,.$$ Then the
inverse map is:
\begin{equation}
\label{ad-stereo1}
\varphi=\ln\Bigl(\varepsilon\frac{x-1}{x+1}\Bigr),\quad
u=-\frac{2y}{(1-x)^{2}}\,.
\end{equation}

If we expand $\mathbb{R}^2$ to the semi-conformal plane $C^2$ by
infinitely distant point and ideal line crossing this point, then
the stereographic map $f$ becomes diffeomorphism $\widetilde{f}$.
This infinitely distant point is the image of point $N$. The ideal
line is the image of straight line belonging to $S^{2}(1)$ and
crossing the pole: $x_{0}=1,\, x_{1}=0$. Do not forget, that we need
to exclude lines $x\neq \pm 1$ from the the plane $\mathbb{R}^2$.

Let us now consider the commutative diagram:
$$
\begin{array}{rcccccl}
S^{2}(1) & & \stackrel{\widetilde{f}}\longrightarrow & & C^{2}\\
& \pi \searrow  & &  \swarrow p \\
 & & \mathbb{R} & & &
\end{array}
$$
The map $p=\pi\circ \widetilde{f}^{-1}:C^{2}\rightarrow R$ is
defined by this diagram. We find the coordinate form of this map:
$$u=-\frac{2y}{(1-x)^{2}}\,.$$
So, the map $p:C^{2}\rightarrow \mathbb{R}$ defines the trivial
principal fibration with the base $\mathbb{R}$ and the structure
group $S_1$.

\textbf{Proposition.} \emph{The map $\widetilde{f}: S^{2}(1)\to
C^{2}$ is conformal.}

\textbf{Proof.} The metric on $G$ induces the metric on $C^{2}$. In
the coordinates $x, y$ it has the form:
\begin{equation}
\label{metric1'}d\widetilde{s}^{2}=-dx^{2}.
\end{equation}

Let us find the metric of semi-Euclidean sphere from the metric on
$C^{2}$. By (\ref{ad-stereo1}) we get $d\varphi=\frac{2}{x^2-1}$.
So, by (\ref{metric1}) and (\ref{metric1'}), we find:
$$
ds^{2}_{1}=\frac{4}{(x^2-1)^2}d\widetilde{s}^{2}.
$$
Hence, the linear element of semi-Euclidean sphere differs from the
linear element of semi-plane by conformal factor and so, the map
$\widetilde{f}$ is conformal. $\Box$

We find the fibers equations on $C^{2}$. The 1-parametric fibers
family of the fibration $(S^{2}(1),\pi,\mathbb{R})$ in the adaptive
coordinates (\ref{adapt1}) is: $u=c,\quad c\in \mathbb{R}.$ By
(\ref{ad-stereo1}) we get the image of this family under map
$\widetilde{f}$:
\begin{equation}\label{fiber1'}
y=-c/2\cdot(x-1)^2.
\end{equation}
The semi-conformal plane is also fibred by this 1-parametric family
of curves. It is parabolas with axis $x=1$ and top $(1,0).$

\section{The projective semi-conformal model of the fibration of sphere $S^2(1)$}

Now we construct the projective semi-conformal model of the sphere
$S^{2}(1)$ and the principal fibration on it. We use normalization
method of A.P.Norden \cite{Nord}. A. P. Shirokov in his work
\cite{Shir} constructed conformal models of Non-Euclidean spaces
with this method.

In a projective space $P_n$ a hypersurface $X_{n-1}$ as an
absolute is called \emph{normalized} if with every point $Q\in X_{n-1}$ is associated:\\
1) line $P_I$ which has the point Q as the only intersection with
the tangent space $T_{n-1}$,\\ 2) linear space $P_{n-2}$ that
belongs
to $T_{n-1}$, but it does not contain the point $Q$.\\
We call them \emph{normals of first and second types}, $P_I$ and
$P_{II}$.

In order that normalization be polar, $P_I$ and $P_{II}$ must be
polar with respect to absolute $X_{n-1}$.

We expand the semi-Euclidean space $_{2}E^{3}_{1}$ to a projective
space $P^{3}$. Here $_{k}E^{n}_{l}$ denotes an $n$-dimensional
semi-Euclidean space with metric tensor of rank $k$, $l$ is the
number of negative inertia index in a quadric form. We consider
homogeneous coordinates $(y_{0}:y_{1}:y_{2}:y_{3})$ in $P^3$, where
$x_{i}=\frac{y_{i}}{y_{3}}\,, i=0,1,2$. Thus $S^{2}(1):\,
x_{0}^{2}-x_{1}^{2}=1$ describes the hyperquadric in $P^{3}$:

\begin{equation}\label{Q}
y_{0}^{2}-y_{1}^{2}-y_{3}^{2}=0\,.
\end{equation}

Here the projective basis $(E_{0},E_{1},E_{2},E_{3})$ is chosen in
the following way. The vertex $E_{0}$ of basis is inside the
hyperquadric. The other vertices $E_{1},E_{2},E_{3}$ are on its
polar plane, $y_{0}=0$. The line $E_0 E_3$ cross the hyperquadric at
poles $N(1:0:0:1)$, $N'(1:0:0:-1)$. Vertices $E_{1},E_{2}$ lie on
the polar of the line $E_0 E_3$. The vertex of the hyperquadric
coincides with the vertex $E_2$.

The stereographic map of projective plane $P^2: y_0=0$ to the
hyperquadric (\ref{Q}) from the pole $N(1:0:0:1)$ is shown on the
picture. Let $U(0:y_{1}:y_{2}:y_{3})\in P^{2}$. If $y_3=0$, then the
line $UN$ belongs to the tangent plane $T_N: y_0-y_3=0$ of the
hyperquadric (\ref{Q}) at the point $N$ and in this case the
intersection point of the line $UN$ with the hyperquadric is not
uniquely determined. If $y_{3}\neq 0$, then the intersection point
of line $UN$ with the hyperquadric is unique. So, we choose the line
$E_{1}E_{2}: y_{3}=0$ as the line at infinity. In the area
$y_{3}\neq 0$ we consider the Cartesian coordinates
$x_{1}=\frac{y_{1}}{y_{3}}, x_{2}=\frac{y_{2}}{y_{3}}\,.$ Then the
plane $\alpha: y_{0}=0, y_{3}\neq 0$ becomes a plane with an affine
structure $A^{2}$. It is possible to introduce the structure of
semi-Euclidean plane $_{1}E^{2}$ with the linear element
\begin{equation}\label{ds0} ds_{0}^{2}=dx_{1}^{2}\,.
\end{equation}
The hyperquadric and the plane $\alpha$ do not intersect or
intersect in two imaginary parallel lines
\begin{equation}\label{q'} x_{1}^{2}=-1\,.
\end{equation}
The restriction of the stereographic projection to the plane
$\alpha$ maps the point $U(0:x_{1}:x_{2}:1)$ into the point $X_1$
\begin{equation}\label{X1}
X_{1}(-1-x_{1}^{2}:2x_{1}:2x_{2}:1-x_{1}^{2})\,.
\end{equation}
So, the Cartesian coordinates $x_i$ can be used as the local
coordinates at the hyperquadric except the point of its intersection
with the tangent plane $T_N$.

\begin{figure}
\begin{center}
\epsfysize=8cm\epsfbox{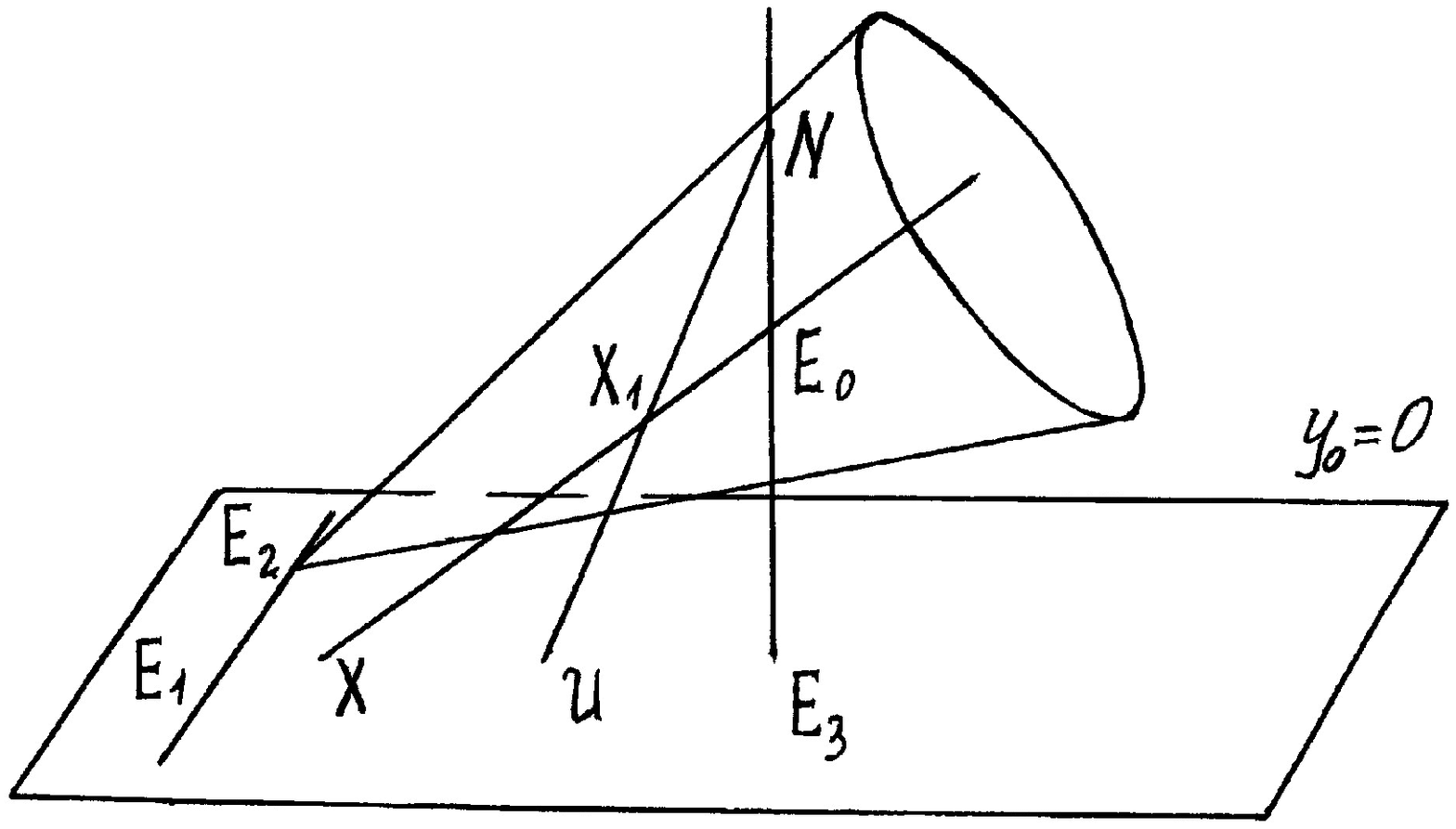}
\end{center}
\end{figure}

We construct an autopolar normalization of the hyperquadric. As a
normal of first type we take lines with fixed center $E_0$ and as a
normal of second type we take their polar lines which belong to the
plane $\alpha$ and cross the vertex $E_2$ of the hyperquadric. The
line $E_{0}X_{1}$ intersects the plane $\alpha$ at the point
$$X(0:2x_{1}:2x_{2}:1-x_{1}^{2})\,.$$ In this normalization the polar
of the point $X$ intersects the plane $\alpha$ on normal $P_{II}$.
Thus for any point $X$ in the plane $\alpha$ there corresponds a
line which does not cross this point. It means that the plane
$\alpha$ is also normalized. Normalization of $\alpha$ is defined by
an absolute quadric (\ref{q'}).

We consider the derivative equations for this normalization. If we
take normals of the first type with fixed center $E_0$, then the
derivative equations (\cite{Nord}, p.204) have the form:
\begin{equation}\label{deriv}
\begin{array}{c}
  \partial_{i}X=Y_{i}+l_{i}X\,, \\
  \nabla_{j}Y_{i}=l_{j}Y_{i}+p_{ji}X\,. \\
\end{array}
\end{equation}
The points $X,\,Y_{i},\,E_{0}$ define a family of projective frames.
Here $Y_{i}$ are generating points of normal $P_{II}$.

We can calculate the values $(X, X),\; (X,Y_i)$ on the plane
$\alpha$ using the quadric form, which is in the left part of
equation (\ref{Q}). So, $(X,X)=-(1+x_{1}^{2})^2.$

Let us find coordinates of metric tensor on plane $\alpha$. So, we
take the Weierstrass standartization

$$(\widetilde{X},\widetilde{X})=-1,\quad
\widetilde{X}=\frac{X}{1+x_{1}^{2}}\,.$$

Then the coordinates of metric tensor are scalar product of partial
derivatives
$g_{ij}=-(\partial_{i}\widetilde{X},\partial_{j}\widetilde{X})$:
$$
(g_{ij})=\left(%
\begin{array}{cc}
  \frac{4}{(1+x_{1}^{2})^2} & 0 \\
  0 & 0 \\
\end{array}\,.%
\right)
$$
So, we get the conformal model of polar normalized plane $\alpha:
y_0=0, y_3\neq0$ with a linear element
\begin{equation}\label{ds}
ds^{2}=\frac{dx_{1}^{2}}{(1+x_{1}^{2})^2}\,.
\end{equation}
It means that this non-Euclidean plane is conformal to
semi-Euclidean plane $_{1}E^{2}$.

The points $X$ and $Y_{i}$ are polarly conjugated, $(X,Y_{i})=0$.
From this equation and the derivative equations (\ref{deriv}) we can
get the non-zero connection coefficients:
$$\Gamma^{1}_{11}=\Gamma^{2}_{12}=\Gamma^{2}_{21}=-\frac{2x_{1}}{1+x_{1}^{2}}\,,\quad
\Gamma^{2}_{11}=\frac{2x_{2}}{1+x_{1}^{2}}\,.$$ The sums
$\Gamma^s_{ks}=\partial_k \ln \frac{c}{(1+x_1^2)^2}$ ($c=const$) are
gradient, so the connection is equiaffine. Curvature tensor has the
following form:
$$R_{121\cdot}^{\quad 2}=-R_{211\cdot}^{\quad 2}=-\frac{4}{(1+x_1^2)^2}\,.$$
Ricci curvature tensor $R_{sk}=R^{\quad i}_{isk\cdot}$ is symmetric:
$R_{11}=\frac{4}{(1+x_1^2)^2}.$ Metric $g_{ij}$ and curvature
$R^{\quad i}_{rsk\cdot}$ tensors are covariantly constant in this
connection: $\nabla_{k}g_{ij}=0,\;\nabla_{l}R^{\quad i}_{rsk\cdot}=
0.$

Geodesic curves in this connection satisfy the following equations:
$$x_2=A(x_1^2-1)+Bx_1 \mbox{ and } x_1=0,$$ where $A$ and $B$ are
arbitrary constants. There are parabolas and lines orthogonal to
axis $Ox_1$.

Let us consider the fibration of this plane by double numbers
subalgebra. We write the equations of fibers of semi-Euclidean
sphere $S^2(1)$ in homogeneous coordinates:
\begin{equation}\label{fiber2}
\left\{%
\begin{array}{ll}
    (y_0-y_1)v-y_{2}=0, \\
    y_{0}^2-y_{1}^2-y_{3}^2=0. \\
\end{array}%
\right.
\end{equation}
This 1-parametric family of curves fibers the hyperquadric and it
defines a fibration on it. The stereographic projection of these
fibers from the pole $N$ to the plane $\alpha$ is:
$$x_2=-v/2\cdot(x_1+1)^2.$$
It is 1-parametric family of parabolas.

\section*{Remark}

We would obtain the similar results for the space of right cosets by
Lie subgroup $H_2$ (it is the subgroup of invertible dual numbers)
and the fibration of the group $G$ by $H_2$. However, $H_{2}$ is a
normal divisor of the group $G$. Therefore, the spaces of right and
left cosets coincide. We have the proposition:

\textbf{Proposition} \emph{The fibration $(G, \pi', M'=G/H_2)$
determined by the formula
\begin{equation}\label{pi2} \pi'(x)=\frac{x_{1}}{x_{0}}
\end{equation}
is a trivial principal fibration over the real line
$\mathbb{R}\backslash\{0\}$. The typical fiber is a plane without a
line. The structure group is $H_2$.}

The equations of fibers are:
\begin{equation}\label{fiber2} ux_{0}-x_{1}=0, \; u\in \mathbb{R}.
\end{equation}
This is a 1-parametric family of planes with common axis $Ox_2$. It
is necessary to remove intersection points of these planes with
planes: $x_{0}\pm x_{1}=0,\;x_0=0$.

The left multiplications $x'=ax$ on the group $G$ preserve fibres if
and only if when $a\in H_2$.

The right multiplications: $x'=xb$ on group $G$ forms the
3-parametric Lie group of linear transformations that preserve the
fibration and induce a 1-parametric group of hyperbolic
transformations on the base with two fixed points $u=\pm1$:
$$
u'=\frac{u+\alpha}{\alpha\cdot u+1},
$$
here $\alpha=\pi(b)$.

The right multiplications include involutions of hyperbolic type:
$$
R(b)=\left(\begin{array}{crc} 0&1&0\\
1&0&0\\
b_{2}&b_{2}&-1\\
\end{array} \right),
 R^{2}={\Id}.
$$
They induce the involution on the base:
$$
u'=\frac{1}{u}.
$$

The similar properties has the fibration of semi-Euclidean sphere by
subgroup $S_2=\{x\in H_2| x_0^2-x_1^2=1\}$ and its models.

\section*{Acknowledgement}
I would like to thank Professor Jiri Vanzura for fruitful
discussions and support with writing this paper. This work was
supported in part by grant No. 201/05/2707 of The Czech Science
Foundation and by the Council of the Czech Government MSM
6198959214.

\end{document}